\documentclass[A4paper,11pt]{article}

\pdfoutput=1

\usepackage{latexsym}
\usepackage{epsfig}
\usepackage{verbatim}
\usepackage{multirow}
\usepackage{fancybox}
\usepackage{shadow}
\usepackage{wrapfig, lipsum}
\usepackage{cite}
\usepackage{amssymb}
\usepackage{amsmath}
\usepackage{bm}
\usepackage{amscd}
\usepackage{graphicx}
\usepackage{xcolor}

\newtheorem{theorem}{\sc Theorem}[section]
\newtheorem{lemma}[theorem]{\sc Lemma}

\newtheorem{proposition}[theorem]{\sc Proposition}

\author{Chu Wang% <-this % stops a space
\thanks{$^{1}$Nokia Bell Labs, 600-700 Mountain Avenue, Murray Hill, New Jersey 07974, 
{\tt chu.wang}@ {\tt nokia.com }}%
}

\title{An Opinion Dynamics Model with Increasing Self-Confidence}

\vspace{.5cm}

\date{}

\begin{document} \maketitle
%\icmlkeywords{Iterated learning, language evolution, opinion dynamics, convergence, sustainable learning}

\vskip 0.3in

%%%%%%%%%%%%%%%%%%%%%%%%%%%%%%%%%%%%%%%%%%%%%%%%%%%%%%%%%%%%%%%%%%%%%%%%%%%%%%%%
\begin{abstract}
We propose an opinion dynamics model in which agents gradually increase their own self-confidence while interacting with each other.
The relations between the newly proposed model and existing works of social learning, inertial opinion dynamics, Bayesian inference, and stochastic multi-armed bandits are demonstrated.
We prove the convergence of the system with the existence of a truth under fixed and periodically changing social networks, 
and obtain tight convergence bounds related to the spectral gap of the graph Laplacian and the maximum total degree centrality, respectively.
In the case of randomly generated social networks, an almost-sure convergence result is obtained.
The dynamics of the model with multiple truths or zero truth is also discussed.
\end{abstract}

\section{Introduction}
In this paper, we propose an opinion dynamics model in which agents move via convex combinations of the positions of neighbors specified by a sequence of graphs.
The key feature of our model is that each agent gradually increases its self-confidence represented by its weight in the convex combination.
The intuition is that by constantly incorporating signals into its own belief, the agent becomes more and more confident about its own opinion.
The motivation and validity of such a scheme are further testified by its relations between related works 
including social learning~\cite{jadbabaie2012non,molavi2015foundations,golub2012homophily,jadbabaie2013information},
Bayesian inference~\cite{box2011bayesian,diaconis1986consistency}, 
multi-armed bandits~\cite{auer2002finite,bubeck2012regret},
as well as opinion dynamics with stubborn or static agents~\cite{hegselmann2015opinion,lorenz2010heterogeneous,chazelle2015}.
We first give a brief introduction of related works, followed by formally presenting the model and its relations to the above areas.
Our main results regarding consensus and convergence rate of the model under fixed and changing social networks are presented in Section \ref{sec:fixed} and Section \ref{sec:changing}, respectively. Section \ref{sec:conclusion} concludes the paper.

Opinion dynamics model and network-based dynamical systems have received a surge of attention lately~\cite{axelrod2006evolution, blondel2005convergence, easley2010networks, castellano2009statistical, blondel2009krause}. 
In these systems, typically, a group of agents will interact by communicating through a sequence of graphs.
The original goal of opinion dynamics is to model the formation and propagation of opinions and knowledge in a crowd of interacting individuals~\cite{hegselmann2002opinion,degroot1974reaching}.
Later the model grows its popularity due to its widespread use in economics and social sciences~\cite{easley2010networks,axelrod2006evolution,chazelle2015}.
Opinion dynamics models are usually diffusive, in the sense that the new opinion comes from a convex combination of the old ones\cite{chazelle2015diffusive}.
Thus, in its matrix form, the dynamics of the opinions is equivalent to repeated stochastic matrix multiplications to the opinions. 
We refer interested readers to~\cite{touri2012product,chazelle2011total,tahbaz2008necessary,cao2008reaching} for more detailed discussions on product of stochastic matrices and averaging processes.

The major question to ask regarding opinion dynamics systems is whether the agents will achieve consensus, a state that all the agents share the same opinion.
In fact, we will show that, with the existence of a static agent (truth), the system achieves consensus under mild assumptions about the network structure.
For fixed graph, convergence is achieved of polynomial order $O(t^{-\nu})$, where $\nu$ is the spectral gap of the graph Laplacian.
For periodically changing graphs, the system converges to consensus asymptotically of order $O(t^{-1/d_\mathrm{max}})$, where $d_\mathrm{max}$ is the maximum total outdegree in a period. 
Both the convergence bounds are tight.
If the graphs are generated randomly, we obtain an almost-sure convergence result.
We note that related opinion dynamics systems usually feature exponential convergence~\cite{blondel2005convergence,cao2008reaching,jadbabaie2012non}. 
The slower convergence in our model originates from the increasing self-confidence: 
larger self-weight after repeated interactions makes the agents reluctant to move. 

\subsection{The Model}
The system consists of $n$ agents represented by scalars or vectors $x_1(t), x_2(t),\dots,x_n(t)$, where $t=0,1,2,\dots$ is the time.
The interactions between agents are captured by a sequence of graphs $G(0), G(1),G(2) \dots$,
where $G(t)$ can be any directed graph over $N$, and self-loop is allowed.
It should be noted that we do not assert constraint of how $G(t)$ is formed, 
and thus $G(t)$ can be fixed, arbitrarily specified, randomly generated, or even coupled with the opinions. 
A truth (static agent) refers to the agent which never interacts with other agents, and thus it is stuck at its initial position eternally.
Though being interesting, it is not necessary for a truth to exist in the system;
if it does, the other mobile agents are usually referred to as the learners.
We use $N_i(t)$ to denote the neighbor set of $i$ at time $t$, and $j\in N_i(j)$ if and only if edge $(i,j)\in G(t)$.
In this case, $i$ gets information from $j$, or equivalently, $j$ influences $i$.
The dynamics of the model is written as
\begin{equation}\label{eq:main}
x_i(t+1)=\frac{w_i(t)x_i(t)+\sum_{j\in N_i(t)}x_j(t)}{w_i(t)+|N_i(t)|},
\end{equation}
where the weight $w_i(t)$ is a scalar associated with agent $i$ at time $t$.
$w_i(t)$ is regarded as $i$'s self-confidence, representing the degree of how much agent $i$ believes in its current opinion.
In spite of various seemingly plausible ways of modeling the dynamics of the self-confidence $w_i(t)$, 
we assume in each step $t$, $w_i(t)$ increases by the number of $i$'s neighbors $|N_i(t)|$:
\begin{equation}\label{eq:weight}
w_i(t+1)=w_i(t)+|N_i(t)|.
\end{equation}
The intuition is that after communicating with $|N_i(t)|$ agents and obtaining $|N_i(t)|$ signals, 
the amount of self-confidence should also increase in that amount.
The modeling of the dynamics of self-confidence \eqref{eq:weight} is further justified by the relations between the proposed model and existing works on social learning, Bayesian inference, inertial opinion dynamics, and multi-armed bandits in Section \ref{sec:relation}, \ref{sec:relation2}, and \ref{sec:relation3}.
In the mean time, out results regarding the proposed model in Section \ref{sec:fixed} and \ref{sec:changing} directly apply to the above fields.

%%%%%%%%%%%%%%%%%%%%%%%%%%%%%%%%%%%%%%%%%%%%%%%%%%%%%%%%%%

\subsection{Relation to Bayesian without Recall Social Learning}\label{sec:relation}
In the framework of social learning, a group of learners (mobile agents) tries to learn the state of the world denoted by a truth (static agent) via a social network.
Recently, Rahimian and Jadbabaie proposed the so-called Bayesian without recall ({\it BWR}) model~\cite{rahimian2016learning}, 
in which the agents are assumed to be rational and memoryless.
In the {\it BWR} model, each agent adopts an initial belief, updates the belief via Bayes's rule based on signals transferred from the truth or other agents, 
while ignoring the mechanism behind the data generating process~\cite{rahimian2015learning, rahimian2015log, rahimian2016naive}.

Assume the belief of agent $i$ is Gaussian distributed $N(\mu_i(t),\sigma_i^2(t))$, 
and the signal being transferred $d_j(t)$ is noisy measurement of $j$'s belief $d_j(t)\sim N(\mu_j(t),\sigma_j^2(t))+\epsilon_j(t)$,
where $\epsilon_j(t)\sim N(0,\sigma^2)$ is independent Gaussian noise.
In \cite{rahimian2016learning}, the explicit update rule is demonstrated as 
\begin{equation}\label{eq:bwr}
\mu_i(t)=\frac{\tau_i(t)\mu_i(t)+\tau\sum_{j\in N_i(t)}\mu_j(t)}{\tau_i(t)+\tau|N_i(t)|},
\end{equation}
where $\tau_i(t):=\sigma_i^{-2}(t)$ and $\tau=\sigma^{-2}$ are inverse variances following update rule
\begin{equation}\label{eq:bwr2}
\tau_i(t+1)=\tau_i(t)+\tau|N_i(t)|,
\end{equation}
By taking expectations on both sides of \eqref{eq:bwr}, it is clear that $\mathbb{E}\mu_i(t)$ together with the weight $\tau_i(t)$ follows the proposed increasing self-confidence model.

%%%%%%%%%%%%%%%%%%%%%%%%%%%%%%%%%%%%%%%%%%%%%%%%%%%%%%%%%%

\subsection{Relation to Inertial Hegselmann-Krause System}\label{sec:relation2}
In the famous Hegselmann-Krause ({\it HK}) system, each agent moves to the mass center of all the agents within a fixed distance $R$~\cite{hegselmann2002opinion}.
Stubborn agent in an {\it HK} system moves toward the mass center of its neighbors by any fraction of length: 
\begin{equation}\label{eq:inertial}
x_i(t+1)=(1-\lambda_i(t))x_i(t)+\frac{\lambda_i(t)}{|N_i(t)|}\sum_{j\in N_i(t)}x_j(t).
\end{equation}
Setting this fraction to zero makes the agent static.
{\it HK} systems with stubborn or static agents attract much attention recently,
for the model is more realistic and addresses issues like symmetry breaking and non-shrinking convex hull~\cite{ghaderi2013opinion,mirtabatabaei2012opinion,chazelle2015}.
The factor $(1-\lambda_i(t))$ in \eqref{eq:inertial} can be regarded as the normalized self-confidence in the proposed model.
In the extreme case of static agents, the corresponding self-confidence is infinity.
We note that in~\cite{chazelle2015}, the inertial $\lambda_i(t)$ 
is endowed with more degrees of freedom and may follow different dynamics other than \eqref{eq:weight}.

The self-confidence $w_i(t)$ should be distinguished from the confidence bound.
The former refers to the self-weight during the convex combination, while the latter is the cut-off threshold of the neighbor set.

%%%%%%%%%%%%%%%%%%%%%%%%%%%%%%%%%%%%%%%%%%%%%%%%%%%%%%%%%%

\subsection{Relation to Bayesian Inference and Multi-armed Bandits}\label{sec:relation3}
To evaluate a quantity $\theta$ by repeated taking noisy measurement $d(t)=\theta+\epsilon(t)$,
one can adopt a Gaussian estimator $N (\mu(t),\sigma^2(t))$ and update it sequentially according to Bayes' rule by formula
\begin{equation}\label{eq:inference}
\mu(t+1)=\frac{\tau(t)\mu(t)+\tau d(t)}{\tau(t)+\tau},
\end{equation}
where the inverse variance $\tau(t)=\sigma^{-2}(t)$ and $\tau=\sigma^{-2}$.
Note that $\mathbb{E}\mu(t)$ forms an increasing self-confidence model of a single truth and a single learner.

The Bayesian multi-armed bandit problem considers $K$ arms with unknown values $\theta_1,\dots,\theta_K$.
A learner sequentially picks arms and obtains noisy feedbacks.
Based on \eqref{eq:inference}, all the $\mathbb{E}\mu_k(t)$ together with the arms form an increasing self-confidence model with $K$ leaners and $K$ truths.
In multi-armed bandits, the goal is either to identify the best arm $\arg\max \mu_k$, or to maximize the expected cumulative rewards $\max \sum_{t\le H} \theta_{k(t)}$,
which lies on whether and how fast $\mathbb{E}\mu_i(t)$ or $\mu_i(t)$ converges to $\theta_k$.
In the language of opinion dynamics, the problem is whether and how fast consensus can be achieved for each pair $\theta_k$ and $\mathbb{E}\mu_k$.

%%%%%%%%%%%%%%%%%%%%%%%%%%%%%%%%%%%%%%%%%%%%%%%%%%%%%%%%%%

\subsection{Dynamics in Matrix Form}
Let $A(t)$ denote the associated adjacency matrix of $G(t)$, thus $a_{ij}(t)=1$ if $(i,j)\in G(t)$, and otherwise $a_{ij}(t)=0$.
We use $D(t)$ to denote the outdegree matrix of $G(t)$: $D(t)$ is a diagonal matrix with its $i$-th diagonal element being the ourdegree of $i$ in $G(t)$.
The weight matrix $W(t):=\mathrm{diag}(w_1(t),w_2(t),\dots,w_n(t))$, can alternatively be defined as $W(t+1)=W(t)+D(t)$ for $t\ge 0$, where $W(0)$ is the initial self-confidence.

We use the letter without the subscript $i$, namely $x(t)$, to denote the column vector $(x_1(t),x_2(t),\dots,x_n(t))^T$.
In this symbol system, the dynamics \eqref{eq:main} is written as
\begin{equation}\label{eq:matrix}
x(t+1)=\left( W(t)+D(t)\right)^{-1}\left( W(t)+A(t)\right)x(t).
\end{equation}
Notice that $D(t)$ is the outdegree matrix of $A(t)$, then $\left( W(t)+D(t)\right)^{-1}( W(t)+A(t))$ is a row-stochastic matrix.
For simplicity, whenever a row $i$ of matrix $W(t)+D(t)$ is zero, the dynamics \eqref{eq:matrix} should be understood in the sense $x_i(t+1)=x_i(t)$,
since a zero row of matrix $W(t)+D(t)$ means the corresponding agent never interacts with anyone.
Without loss of generality, we assume the initial self-confidence is 0 for each agent throughout this paper.
For any matrix $M$ used in the paper, we use the small letter with double subscripts $m_{ij}$ to denote the $(i,j)$-th element of $M$.

%%%%%%%%%%%%%%%%%%%%%%%%%%%%%%%%%%%%%%%%%%%%%%%%%%%%%%%%%%
\section{Truth Seeking in Fixed Social Network}\label{sec:fixed}
In this section, we consider the increasing self-confidence model with a truth and fixed graph $G(t)=G$.
In this case, $A(t)=A$, $D(t)=D$, $W(t)=tD$, 
and the matrix-form dynamics \eqref{eq:matrix} becomes
\begin{equation}\label{eq:fixed_update}
x(t+1)=((t+1)D)^{-1}(tD+A)x(t).
\end{equation}
Notice that, unlike the DeGroot model~\cite{degroot1974reaching}, 
a sequence of the same graph does not lead to repeated multiplication of the same stochastic matrix.
This is because the self-confidence of an agent will increase after it receives signals, 
and thus the update rule \eqref{eq:main} does not remain the same for different $t$.

Without loss of generality, we assume that agent 1 is the truth that stays at the origin forever. 
Thus, $a_{1,j}=0$ for all $j>1$. 
Notice that the value of $a_{11}$ will not affect the dynamics, hence we assume $a_{11}=1$ for convenience.
By recursively adopting \eqref{eq:fixed_update}, we have
\begin{eqnarray}\label{eq:fixed_recursive}
&&(t+1)Dx(t+1)=Ax(t)+tDx(t)\nonumber\\
&=&Ax(t)+Ax(t-1)+(t-1)Dx(t-1)\nonumber\\
&=&Ax(t)+\cdots+Ax(0).
\end{eqnarray}
Let $S(t):=x(t)+x(t-1)+\cdots+x(0)$, then from \eqref{eq:fixed_recursive}, the dynamics of $S(t)$ is written as
\begin{equation}\label{eq:fixed_S}
S(t+1)=\left(I+\frac{D^{-1}A}{t+1}\right)S(t).
\end{equation}
By repeatedly adopting \eqref{eq:fixed_S}, we obtain a formula of $x(t)$ as
\begin{eqnarray}\label{eq:product}
x(t)&=&S(t)-S(t-1)=(tD)^{-1}AS(t-1)\nonumber\\
&=&\frac{D^{-1}A}{t}\prod_{s=1}^{t-1}\left(I+\frac{D^{-1}A}{s}\right)x(0).
\end{eqnarray}
Intuitively, if we replace the matrix $D^{-1}A$ with a real number $0<\rho<1$, then the product of matrices in \eqref{eq:product} becomes 
\begin{equation}\label{eq:rho}
\frac{\rho}{t}\prod_{s=1}^{t-1}\left(1+\frac{\rho}{s}\right)\le \frac{\rho}{t}\exp\left(\rho\sum_{s=1}^{t-1}\frac{1}{s}\right)=O(t^{\rho-1}).
\end{equation}
Therefore the vanishing speed is of polynomial order $O(t^{\rho-1})$.
Back to the dynamics of \eqref{eq:product}, we claim that whether $x_i(t)$ vanishes (converges to the truth) depends on the difference between 1 and the modulus of the second largest eigenvalue of $D^{-1}A$, which in fact is the spectral gap of $G$. 
Formally, we claim:
\medskip
\begin{theorem}\label{th:fixed}
If the graph of the system $G_t=G$ is fixed, and each leaner has a path to the truth in $G$, 
then the system converges to the truth in polynomial order $O(t^{-\nu})$, where $\nu>0$ is the spectral gap of $G$.
\end{theorem}
\medskip
Notice that under the assumption that each learner has a path to the truth, the outdegree of each agent is positive, and thus $D$ is invertible.

\begin{proof}
The proof proceeds in two stages. We will first show that $\nu>0$, then we will prove the convergence and estimate the convergence rate.
Since $D^{-1}A$ is a row-stochastic matrix, therefore 1 is its largest eigenvalue. 
To prove that $\nu>0$, it is sufficient to prove that $D^{-1}A$ does not have other eigenvalues with modulus 1.
By regarding the truth and the learners as two groups, it is clear that $D^{-1}A$ is a block lower-triangular matrix. 
We use $E$ and $B$ to denote the remaining matrices by removing the first row and the first column of $D$ and $A$, respectively.
What is left to show is that $E^{-1}B$ does not have an eigenvalue $\lambda^*$ on the unit circle in the complex plane.

Suppose otherwise $E^{-1}Bu=\lambda^* u$ for a non-zero vector $u$, then $(\lambda^* E-B)u=0$. 
Recall that $D$ is the outdegree matrix of $A$, we have
\begin{equation}\label{eq:tau}
|\lambda^* e_{ii}-b_{ii}|\ge e_{ii}-b_{ii}\ge \sum_{j\neq i} b_{ij},
\end{equation}
and the second inequality is strict if $i$ has an edge to the truth.
Suppose $|u_k|=\max_i|u_i|$, then from $(\lambda^* E-B)u=0$ we have
\begin{eqnarray}\label{eq:equal}
0&=&\big{|}(\lambda^* e_{ii}-b_{ii})u_i-\sum_{j\neq i}b_{ij}u_j\big{|} \nonumber\\
&\ge& |\lambda^* e_{ii}-b_{ii}||u_i|-\sum_{j\neq i}|b_{ij}||u_j|\ge0.
\end{eqnarray}
One conclusion from \eqref{eq:equal} is that equalities hold in \eqref{eq:tau}, which means $i$ does not have an edge to the truth.
Furthermore, \eqref{eq:equal} implies that for all $b_{ij}\neq 0$, $|u_j|=|u_i|$ also has the largest modulus.
Therefore, by repeating the same argument of $i$ to $j$, it is clear that $j$, and hence any learner reachable from $i$ in $G$, can not have edge to the truth.
This contradicts the assumption in Theorem \ref{th:fixed} that each leaner has a path to the truth.
Therefore, $|\lambda^*|<1$ and thus we have proved that $\nu>0$.

For the convergence of the system, let $\tilde{x}(t)=(x_2(t),\dots,x_n(t))^T$ denote the dynamics of the learners, then it follows the update rule:
\begin{equation}\label{eq:fixed_learner}
\tilde{x}(t)=\frac{E^{-1}B}{t}\prod_{s=1}^{t-1}\left(I+\frac{E^{-1}B}{s}\right)\tilde{x}(0).
\end{equation}
Now let $U$ be an invertible matrix such that the similarity transformation $U^{-1}E^{-1}BU$ is the Jordan normal form of $E^{-1}B$,
then for the vanishing speed of the product of matrices in \eqref{eq:fixed_learner}, it is equivalent to estimate the vanishing speed of
\begin{equation*}
Q^{(\Lambda)}(t):=\frac{\Lambda}{t}\left(I_{k\times k}+\frac{\Lambda}{t-1}\right)\cdots\left(I_{k\times k}+\Lambda\right),
\end{equation*}
where $\Lambda$ is any Jordan block of matrix $U^{-1}E^{-1}BU$:
\begin{equation}\label{eq:Q_lambda}
\Lambda=\lambda I_{k\times k}+J,~~ J=\begin{pmatrix}
\bm{0}_{(k-1)\times 1} & I_{(k-1)\times(k-1)} \\
0 & \bm{0}_{1 \times (k-1)}
\end{pmatrix},
\end{equation}
and $\lambda$ is an eigenvalue of $E^{-1}B$.
Since $J^{k-1}=\bm{0}_{k\times k}$, and it commutes with multiples of $I_{k\times k}$ in the definition of $Q^{(\Lambda)}(t)$ in \eqref{eq:Q_lambda},
then for any $\lambda\neq 0$,
$$\max_{ij}|q^{(\Lambda)}_{ij}(t)|=O\left(\frac{|\lambda|}{t}\prod_{s=1}^{t-1}\left(1+\frac{|\lambda|}{s}\right)\right)=O(t^{|\lambda|-1}).$$
Notice that $U$ only depends on the graph $G$.
As a result, 
$$\max_{i}|x_i(t)|=\max_{\lambda} O(t^{|\lambda|-1})=O(t^{\max|\lambda|-1})=O(t^{-\nu}),$$
which completes the proof.
\end{proof}
A prevailing assumption of the graph structure in social learning is the network being strongly connected~\cite{jadbabaie2012non,jadbabaie2013information},
which guarantees the truth to be reachable by each learner. 
We note that the reachability of truth is similar to the definition of rooted graph in~\cite{cao2008reaching},
in which Cao {\it et al.} carefully analyzed the convergence of the system with fixed self-confidence.

In addition to the convergence result, we prove that the bound of convergence is tight by constructing the following system.
The initial positions $\tilde{x}_i(0)=1$ for $i\ge 2$;
the outdegree of each learner is the same number $d\ge 1$;
each leaner has one edge to the truth;
and $b_{ij}=1$ if and only if $j=i,i+1,\dots,i+d-2$, where the subscript is understood modulo $(n-1)$.
In the extreme case $d=1$, $B$ is a zero matrix.
Under this construction, $E=dI_{(n-1)\times(n-1)}$, and $E^{-1}B$ is a circulant matrix whose eigenvalues are straightforward to get: 
$\lambda_k=(1+\omega_j+\dots+\omega_j^{d-2})/d$, where $\omega_j=\exp(2\pi \bm{i}j/(n-1))$ is the $(n-1)$-th root of unity and $\bm{i}$ stands for the imaginary unit.
It is clear to check that $|\lambda_k|\le 1-1/d$ and thus $\nu=1/d$.

On the other hand, since each learner starts at the same position and talks to the same number of learners, we have $x_2(t)=\dots=x_n(t)$.
In view of \eqref{eq:product},
$$
x_{2}(t) =\prod_{s=1}^{t}\left(1-\frac{1}{sd}\right)=\exp\left(\sum_{s=1}^t \log \left(1-\frac{1}{sd}\right)\right).
$$
From Taylor expansion, we have $\log(1-c)\ge -c-c^2/2$ for any $0\le c<1$. 
Therefore
$$
x_{2}(t)\ge \exp\left(-\frac{1}{d}\sum_{s=1}^t \frac{1}{s}-\frac{1}{2d^2}\sum_{s=1}^t \frac{1}{s^2}\right)=\Omega (t^{-\frac{1}{d}})=\Omega (t^{-\nu}).
$$
We have proved
\medskip
\begin{proposition}\label{prop}
For any $n$, there exists a graph $G$ with spectral gap $\nu$ such that $\|x(t)\|_\infty=\Omega(t^{-\nu})$.
\end{proposition}
\medskip

%%%%%%%%%%%%%%%%%%%%%%%%%%%%%%%%%%%%%%%%%%%%%%%%%%%%%%%%%%
\section{Changing Social Netwotks}\label{sec:changing}
In this section, we consider the increasing self-confidence model in changing social networks.
In order to achieve consensus, information from the truth should be well spread to the learners.
Indeed, if an agent is not able to get signals from the truth constantly, its dynamics should eventually be free from the influence of the truth.
This intuition leads to the definition and analysis of {\it influence indicator}.
%%%%%%%%%%%%%%%%%%%%%%%%%%%%%%%%%%%%%%%%%%%%%%%%%%%%%%%%%%

\subsection{The Influence Indicator}
We use $P(t)$ to denote the remaining matrix by removing the first row and the first column of $D^{-1}(t)A(t)$.
Recall that whenever $d_{ii}(t)=0$, $p_{ij}(t)$ is set to 0 for all $j$.
Note that $P(t)$ is a sub-stochastic matrix since its row-sums are no greater than 1.
If at time $t$, agent $i$ has an edge pointing to the truth, then $d_{ii}(t)=1+a_{12}(t)+a_{13}(t)+\cdots+a_{1n}(t)$.
Thus, the corresponding row-sum is strictly less than one, implying non-zero impact from the truth.
Indeed, if the row-sum is exactly 1, then the dynamics of the corresponding agent is completely determined by only the learners.
Therefore, the difference between 1 and each row-sum of $P(t)$ is an indicator of the influence from the truth.
Formally, we define the influence indicator from time $s$ to time $t$, denoted by $\alpha(t:s)$, as:
\begin{equation}\label{eq:leak}
\alpha(t:s)=\xi-P(t:s)\xi,
\end{equation}
where $\xi$ is the all-one column vector.
When $t=s+1$, \eqref{eq:leak} is reduced to 
\begin{equation}\label{eq:leak_single}
\alpha(s)=\xi-P(s)\xi.
\end{equation}
Notice that $\alpha(t)=0$ indicates $1\notin N_i(t)$.
By repeatedly left-multiplying $P(s+1),P(s+2),\dots,P(t-1)$ to both sides of \eqref{eq:leak_single}, 
we obtain
\begin{equation}\label{eq:leak_chain}
\alpha(t:s)=\sum_{k=s}^{t-1} P(t:k)\alpha(k),
\end{equation}
which builds up the relation between the single-step indicator $\alpha(k)$ and the multi-step indicator $\alpha(t:s)$.

The dynamics of the learners $\tilde{x}(t)$ is determined by the product of $P(s)$ for $s\le t$: $\tilde{x}(t)=P(t:0)\tilde{x}(0)$, and thus
\begin{equation}\label{eq:x_t}
\|\tilde{x}(t)\|_\infty\le \|P(t:0)\|_\infty\|\tilde{x}(0)\|_\infty.
\end{equation}
On the other hand, taking infinity norm on both sides of \eqref{eq:leak} yields
\begin{equation}\label{eq:norm}
\|P(t:s)\|_\infty= 1-\min_i \alpha_i(t:s).
\end{equation}
Therefore, whether $\tilde{x}$ converges to the truth depends on the minimum value of the indicator.

To estimate the indicator, notice that matrix $P(t)$ is non-negative, thus $p_{ii}(t:k)\ge p_{ii}(t-1)\dots p_{ii}(k)$.
In addition, for any learner $i$ and time $r$, $p_{ii}(r)=(w_i(r)+a_{i1}(r))/w_i(t+1)\ge w_i(r)/w_i(r+1)$.
Therefore, for any $k\ge s$,
\begin{equation} \label{eq:diagonal}
p_{ii}(t:k) \ge \prod_{r=k}^{t-1} \frac{w_{ii}(r)}{w_{ii}(r+1)}=\frac{w_{ii}(k)}{w_{ii}(t)}\ge \frac{w_{ii}(s)}{w_{ii}(t)}.
\end{equation}
Again, when $w_{ii}(t)=0$, then $w_{ii}(s)=0$, and \eqref{eq:diagonal} should be understood as the trivial inequality $p_{ii}(t:k)\ge 0$.
In view of \eqref{eq:leak_chain} , \eqref{eq:diagonal}, and notice that the matrices and vectors involved are all non-negative,
we obtain a lower bound estimate for the indicator. 
Formally, we have
\medskip
\begin{lemma}\label{lemma}
For any $t>s\ge 0$, the following inequality of the influence indicator holds:
\begin{equation}\label{leak_bound}
\alpha(t:s)\ge W^{-1}(t)W(s) \sum_{k=s}^{t-1}\alpha(k).
\end{equation}
\end{lemma}
\medskip
The inequality \eqref{leak_bound} is element-wise. 
We will also use matrix inequality in the same sense for the rest of this paper.

\iffalse
\begin{eqnarray} \label{eq:diagonal}
&&p_{ii}(k:s) \ge  \prod_{r=s}^{k-1} p_{ii}(r)=\prod_{r=s}^{k-1} \frac{w_{ii}(r)+a_{ii}(r)}{w_{ii}(r+1)}\nonumber\\
&=&\prod_{r=s}^{k-1} \frac{w_{ii}(r)}{w_{ii}(r+1)}=\frac{w_{ii}(s)}{w_{ii}(k)}\ge \frac{w_{ii}(s)}{w_{ii}(t)},
\end{eqnarray}
\fi

%%%%%%%%%%%%%%%%%%%%%%%%%%%%%%%%%%%%%%%%%%%%%%%%%%%%%%%%%%

\subsection{Periodic Graph Sequence}
In this subsection, we consider the case when the graph sequence is periodic: $G_{t+T}=G_{t}$ for $t\ge 0$, where $T\ge 1$ is the period.
Note that the special case $T=1$ reduces to the fixed graph scenario.
We define the total outdegree of agent $i$ in a period $d_{i}:=d_{ii}(1)+\dots+d_{ii}(T)$, and its maximum $d_{\mathrm{max}}:=\min_{2\le i\le n} d_i$.
Since the graph sequence is periodic, the self-confidence $w_{ii}(kT)=kd_i$ grows linearly.
We first state our main result in this subsection:
\medskip
\begin{theorem}\label{th:periodic}
If the graph sequence $G_t$ is periodic, and each learner has at least one edge to the truth in a period,
then the system will converge to the truth in the order $O(t^{-1/d_{\mathrm{max}}})$, and the bound is tight.
\end{theorem}
\medskip
\begin{proof}
Since each learner has at least one edge to the truth in a period, then $a_{i1}(kT)+a_{i1}(kT+1)+\dots+a_{i1}((k+1)T-1)\ge 1$ for any $2\le i\le n$,
and thus each element of the indicator $\alpha((k+1)T:kT)$ is positive.
More precisely,  
$$\alpha_i(kT+l)=\frac{a_{i1}(kT+l)}{w_{ii}(kT+l)}\ge \frac{a_{i1}(kT+l)}{w_{ii}((k+1)T)}=\frac{a_{i1}(kT+l)}{(k+1)d_{i}},$$
and thus
$$\sum_{r=kT}^{(k+1)T-1}\alpha(r) \ge \frac{1}{(k+1)d_i}\sum_{r=kT}^{(k+1)T-1}a_{i1}(r)\ge \frac{1}{(k+1)d_i}.$$
Based on Lemma \ref{lemma},  for $k\ge 1$:
$$
\min_{2\le i\le n}\alpha_i((k+1)T:kT)\ge \frac{k}{(k+1)^2d_{\mathrm{max}}}\ge\frac{1}{(k+3)d_{\mathrm{max}}}.
$$
In view of the last inequality and \eqref{eq:norm},
\begin{equation}\label{eq:Bst}
\|P((k+1)T:kT)\|_\infty\le1-\frac{1}{(k+3)d_{\mathrm{max}}}.
\end{equation}

Now let $t=rT$, by the sub-multiplicativity of the matrix infinity norm,
\begin{eqnarray}
&&\|P(rT:0)\|_\infty\le \prod_{k=1}^{r-1}\|P((k+1)T:kT)\|_\infty \nonumber\\ 
&\le&\prod_{k=4}^{r+2}\left(1-\frac{1}{kd_{\mathrm{max}}}\right)\le 
\exp\left(-\frac{1}{d_{\mathrm{max}}}\sum_{k=4}^{r+2}\frac{1}{k}\right)\nonumber\\
&\le& O(r^{-\frac{1}{d_{\mathrm{max}}}})= O(t^{-\frac{1}{d_{\mathrm{max}}}}).
\end{eqnarray}
This completes the proof of convergence.

For any positive integer $d_{\mathrm{max}}$ and period $T$, 
we build a system with graph sequence $G(t)$ such that $G(kT)$ is the graph in the proof of Proposition \ref{prop}, 
and $G(kT+r)$ is an empty graph for $1\le r\le T-1$.
Then we have 
$$\|x(t)\|_\infty=\Omega((t/T)^{-\frac{1}{d_{\mathrm{max}}}})=\Omega(t^{-\frac{1}{d_{\mathrm{max}}}}),$$
which completes the proof.
\end{proof}
Note that $d_{\mathrm{max}}$ can be interpreted as the maximum total degree centrality of all the learners,
which indicates a slower convergence rate when highly-important leaner exists.
%%%%%%%%%%%%%%%%%%%%%%%%%%%%%%%%%%%%%%%%%%%%%%%%%%%%%%%%%%%%%%%%%%%%%%%%%%%%%%%%%%%%%%%%%%%%%%%%%%%%%%%%%%%%%%%%%%%%%%%%%%%%%%%%%%%%%%%%%%%%%%%%%%%%%%%%%%%%%%%%%%%%%%%%%%%%%

\subsection{Random Graph Sequence}
Now we consider the scenario when the graph $G(t)$ is randomly sampled.
There are various ways of generating a random graph, and we adopt the following scheme.
In each step, every learner randomly picks $d_i$ agents as its neighbors, where $1\le d_i\le n$ is a fixed integer, and we define $d_{\mathrm{max}}=\max_{2\le i\le n}d_i$.
We will show an almost-sure convergence result:
\medskip
\begin{theorem}\label{th:random}
Assume in each step, each learner $i$ independently picks $d_i$ agents as its neighbors uniformly at random, then almost surely, the system converges to the truth,
and the convergence rate is polynomial in $t$.
\end{theorem}
\medskip
\begin{proof}
Since in each step, the outdegree of agent $i$ is the fixed number $d_i$, then the self-confidence $w_{ii}(t)=d_it$.
In view of \eqref{eq:norm}, 
\begin{equation}\label{eq:random_infinity}
\|P(t)\|_\infty=1-\min_i \alpha_i(t)\le1-\frac{1}{d_\mathrm{max}t}\min_i a_{i1}(t).
\end{equation}
Notice that $a_{i1}(t)$ is a Bernoulli random variable with $\mathbb{P}[a_{i1}(t)=1]=d_i/n$.
Define random process $\beta(t)=\min_i a_{i1}(t)$.
Since $a_{i1}(t)$ and $a_{j1}(t)$ are independent for $i\neq j$, then $\beta(t)$ is again a Bernoulli random variable, and
$$q:\mathbb{P}[\beta(t)=1]=d_2d_3\dots d_n/n^{n-1}.$$
Based on \eqref{eq:random_infinity},
$$\hspace*{-0.2cm}
\|B(t:0)\|_\infty \le \prod_{s=1}^{t}\left(1-\frac{\beta(s)}{sd_\mathrm{max}}\right)\le \exp\left(-\frac{1}{d_{\mathrm{max}}}\sum_{s=1}^t\frac{\beta(s)}{s}\right).
$$
We only need to show that, with probability one, the sum of the series $\beta(s)/s$ goes to infinity.

Let $\gamma(s):=\beta(s)-\mathbb{E}\beta(s)$, then $\mathbb{E}\gamma(s)=0$ and $\mathrm{var}~[\gamma(s)]=\mathrm{var}~[\beta(s)]=q(1-q)$.
Define the random process 
$$V(s):=\sum_{k=1}^s \gamma(s)/s$$ and $\mathcal{F}(s)$ the sigma algebra generated by $V(0),\dots,V(s)$.
Then
$$\mathbb{E}[V_{s+1}|\mathcal{F}(s)]=\frac{\mathbb{E}v_s}{s}=0,$$
hence $V(s)$ is a martingale.
In addition, 
$$\mathbb{E}V_s^2=\sum_{k=1}^s \frac{\mathbb{E}v^2(s)}{s^2}\le \sum_{k=1}^s \frac{1}{s^2}<\infty.$$
Therefore by Martingale convergence theorem,  the random variable $V_{\infty}=:\lim_{s\rightarrow\infty}V_s$ exists and has finite variance~\cite{williams1991probability}.
Now let $\beta:= \mathbb{E}\beta(s)$, then almost surely
\begin{eqnarray*}
\|B(t:0)\|_\infty &\le&  \exp\left(-\frac{1}{d_{\mathrm{max}}}\sum_{s=1}^t\frac{\beta+\gamma(s)}{s}\right)\\
&\le& O(t^{-\frac{\beta}{d_{\mathrm{max}}}}),
\end{eqnarray*}
which completes the proof.
\end{proof}

\begin{figure}[htp!]
\centering
\includegraphics[width=0.7\textwidth]{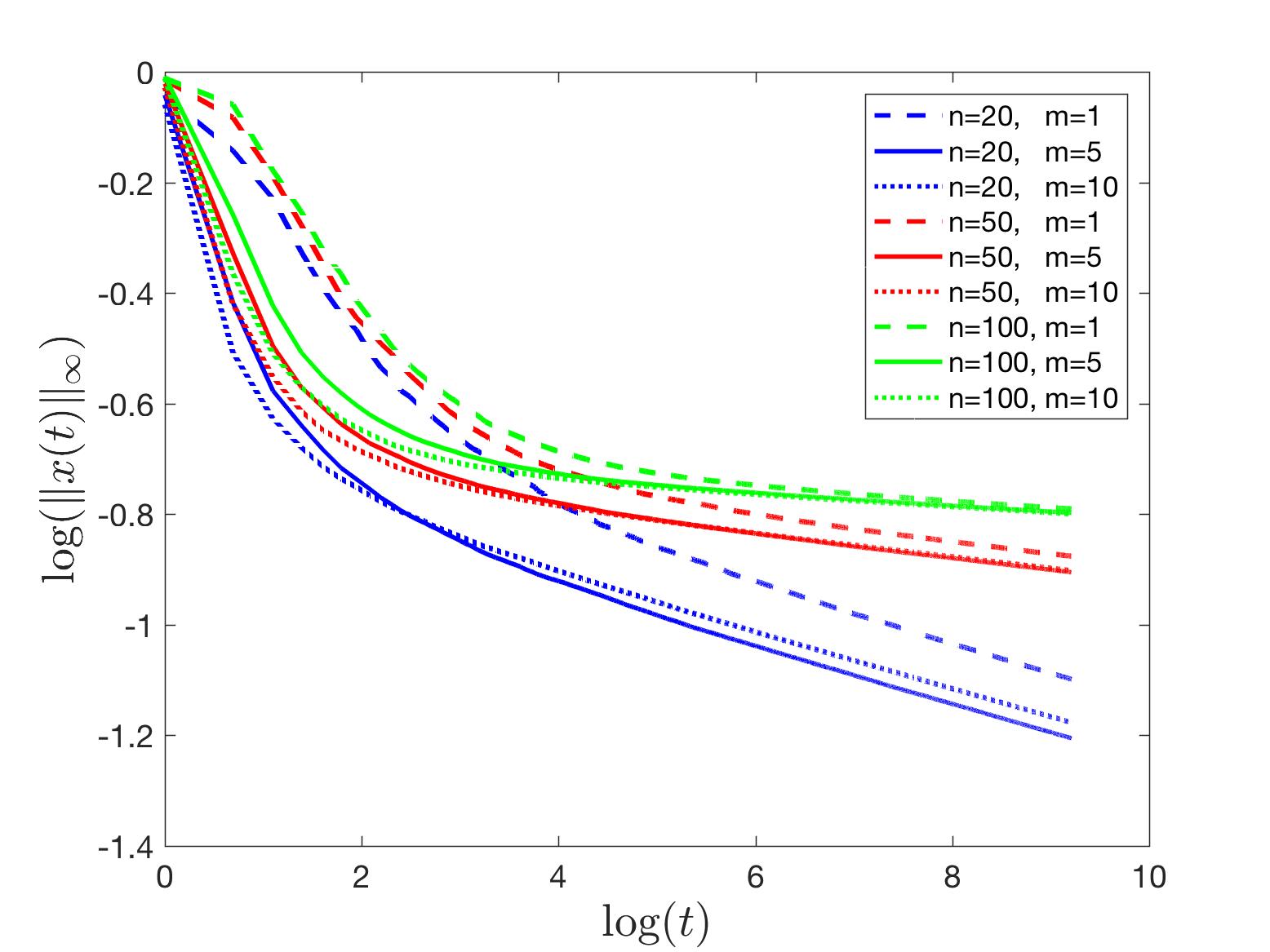}
\caption{Log-log curve of $\|x(t)\|_\infty$ and $t$ for $n=20,50,100$ and $m=1,5,10$. The curves eventually become straight lines with negative slopes.\label{fig}}
\end{figure}

We simulate the system with $n=20,50,100$ agents with $m=1,5,10$.
The initial positions of the learners are uniformly sampled in the unit interval.
For each pair $(n,m)$, we simulate 100 independent systems and calculate the averaged error $\|x(t)\|_\infty$.
The log-log curve of each case is demonstrated in Figure \ref{fig}.
It is clear that the curves eventually become straight lines, indicating a convergence rate of polynomial order.

We further use linear regression to get the slope $f_{n,m}$ of each curve after it becomes steady,
and obtain a rough relation: $f_{n,m}\approx-1/n$, which does not depend on $m$.
For learner $i$, when $m$ increases, it is more likely for the truth to be a neighbor of $i$, which contributes the convergence. 
But in the mean time, there are more learners in the neighbor set of $i$, which harms the convergence since the information from other learners is less perfect compared to the information from the truth.

\subsection{Discussions: Multiple Truths and Zero Truth}
The presence of multiple truths will immediately complicates the behavior of the system.
For example, if a learner communicates with truth 1 enough times in order to be in the vicinity of truth 1, 
and then starts to communicate with truth 2 and does the same,
by repeating this process, the learner will oscillate between the two truths forever.
If the system contains no truth, then the previous example could still happen: we only need to replace each truth by two colliding agents.
Nevertheless, from the proof of Theorem \ref{th:fixed}, it is clear that the components of each learner that is perpendicular to $\mathcal{S}$ should vanish in polynomial order, where $\mathcal{S}$ is the space spanned by all the truths.
The dynamics in the perpendicular space $\mathcal{S}$ is more involved in the graph sequence $G(t)$.
We note that more powerful techniques are required for such general cases.

\section{Conclusion}\label{sec:conclusion}
In this paper, we proposed an opinion dynamics model with increasing self-confidence.
The growing confidence of an agent after it repeatedly communicates with others is reflected in its increasing self-weight.
We proved that, with fixed or periodically changing social network and a single truth, 
the system achieves consensus asymptotically and a tight convergence rate of polynomial order is obtained. 
If each learner randomly selects a fixed number of neighbors, then the system is proved to converge to the truth almost surely.
We also discussed the behavior of the model when zero truth or multiple truths are present, 
which requires delicate analysis in the future.

\footnotesize{
\bibliography{refer}

\begin{thebibliography}{10}

\bibitem{jadbabaie2012non}
Ali Jadbabaie, Pooya Molavi, Alvaro Sandroni, and Alireza Tahbaz-Salehi.
\newblock Non-bayesian social learning.
\newblock {\em Games and Economic Behavior}, 76(1):210--225, 2012.

\bibitem{molavi2015foundations}
Pooya Molavi, Alireza Tahbaz-Salehi, and Ali Jadbabaie.
\newblock Foundations of non-bayesian social learning.
\newblock {\em Columbia Business School Research Paper}, 2015.

\bibitem{golub2012homophily}
Benjamin Golub and Matthew~O Jackson.
\newblock How homophily affects the speed of learning and best response
  dynamics.
\newblock 2012.

\bibitem{jadbabaie2013information}
Ali Jadbabaie, Pooya Molavi, and Alireza Tahbaz-Salehi.
\newblock Information heterogeneity and the speed of learning in social
  networks.
\newblock {\em Columbia Business School Research Paper}, (13-28), 2013.

\bibitem{box2011bayesian}
George~EP Box and George~C Tiao.
\newblock {\em Bayesian inference in statistical analysis}, volume~40.
\newblock John Wiley \& Sons, 2011.

\bibitem{diaconis1986consistency}
Persi Diaconis and David Freedman.
\newblock On the consistency of bayes estimates.
\newblock {\em The Annals of Statistics}, pages 1--26, 1986.

\bibitem{auer2002finite}
Peter Auer, Nicolo Cesa-Bianchi, and Paul Fischer.
\newblock Finite-time analysis of the multiarmed bandit problem.
\newblock {\em Machine learning}, 47(2-3):235--256, 2002.

\bibitem{bubeck2012regret}
S{\'e}bastien Bubeck and Nicolo Cesa-Bianchi.
\newblock Regret analysis of stochastic and nonstochastic multi-armed bandit
  problems.
\newblock {\em arXiv preprint arXiv:1204.5721}, 2012.

\bibitem{hegselmann2015opinion}
Rainer Hegselmann and Ulrich Krause.
\newblock Opinion dynamics under the influence of radical groups, charismatic
  leaders, and other constant signals: A simple unifying model.
\newblock {\em Networks and Heterogeneous Media}, 10(3):477--509, 2015.

\bibitem{lorenz2010heterogeneous}
Jan Lorenz.
\newblock Heterogeneous bounds of confidence: meet, discuss and find consensus!
\newblock {\em Complexity}, 15(4):43--52, 2010.

\bibitem{chazelle2015}
Bernard Chazelle and Chu Wang.
\newblock Inertial {H}egselmann-{K}rause systems.
\newblock In {\em Proceedings of the IEEE American Control Conference (ACC)},
  pages 1936--1941, 2016.

\bibitem{axelrod2006evolution}
Robert~M Axelrod.
\newblock {\em The evolution of cooperation}.
\newblock Basic books, 2006.

\bibitem{blondel2005convergence}
Vincent Blondel, Julien~M Hendrickx, Alex Olshevsky, J~Tsitsiklis, et~al.
\newblock Convergence in multiagent coordination, consensus, and flocking.
\newblock In {\em IEEE Conference on Decision and Control}, volume~44, page
  2996. IEEE; 1998, 2005.

\bibitem{easley2010networks}
David Easley and Jon Kleinberg.
\newblock {\em Networks, crowds, and markets: Reasoning about a highly
  connected world}.
\newblock Cambridge University Press, 2010.

\bibitem{castellano2009statistical}
Claudio Castellano, Santo Fortunato, and Vittorio Loreto.
\newblock Statistical physics of social dynamics.
\newblock {\em Reviews of modern physics}, 81(2):591, 2009.

\bibitem{blondel2009krause}
Vincent~D Blondel, Julien~M Hendrickx, and John~N Tsitsiklis.
\newblock On krause's multi-agent consensus model with state-dependent
  connectivity.
\newblock {\em IEEE transactions on Automatic Control}, 54(11):2586--2597,
  2009.

\bibitem{hegselmann2002opinion}
Rainer Hegselmann and Ulrich Krause.
\newblock Opinion dynamics and bounded confidence models, analysis, and
  simulation.
\newblock {\em Journal of Artificial Societies and Social Simulation}, 5(3),
  2002.

\bibitem{degroot1974reaching}
Morris~H DeGroot.
\newblock Reaching a consensus.
\newblock {\em Journal of the American Statistical Association},
  69(345):118--121, 1974.

\bibitem{chazelle2015diffusive}
Bernard Chazelle.
\newblock Diffusive influence systems.
\newblock {\em SIAM Journal on Computing}, 44(5):1403--1442, 2015.

\bibitem{touri2012product}
Behrouz Touri.
\newblock {\em Product of random stochastic matrices and distributed
  averaging}.
\newblock Springer Science \& Business Media, 2012.

\bibitem{chazelle2011total}
Bernard Chazelle.
\newblock The total s-energy of a multiagent system.
\newblock {\em SIAM Journal on Control and Optimization}, 49(4):1680--1706,
  2011.

\bibitem{tahbaz2008necessary}
Alireza Tahbaz-Salehi and Ali Jadbabaie.
\newblock A necessary and sufficient condition for consensus over random
  networks.
\newblock {\em IEEE Transactions on Automatic Control}, 53(3):791--795, 2008.

\bibitem{cao2008reaching}
Ming Cao, A~Stephen Morse, and Brian~DO Anderson.
\newblock Reaching a consensus in a dynamically changing environment: A
  graphical approach.
\newblock {\em SIAM Journal on Control and Optimization}, 47(2):575--600, 2008.

\bibitem{rahimian2016learning}
M~Amin Rahimian and Ali Jadbabaie.
\newblock Learning without recall from actions of neighbors.
\newblock In {\em 2016 American Control Conference (ACC)}, pages 1060--1065.
  IEEE, 2016.

\bibitem{rahimian2015learning}
Mohammad~Amin Rahimian, Shahin Shahrampour, and Ali Jadbabaie.
\newblock Learning without recall by random walks on directed graphs.
\newblock In {\em 2015 54th IEEE Conference on Decision and Control (CDC)},
  pages 5538--5543. IEEE, 2015.

\bibitem{rahimian2015log}
Mohammad~Amin Rahimian et~al.
\newblock Learning without recall: A case for log-linear learning.
\newblock {\em IFAC-PapersOnLine}, 48(22):46--51, 2015.

\bibitem{rahimian2016naive}
Mohammad~Amin Rahimian and Ali Jadbabaie.
\newblock Naive social learning in ising networks.
\newblock In {\em 2016 American Control Conference (ACC)}, pages 1088--1093.
  IEEE, 2016.

\bibitem{ghaderi2013opinion}
Javad Ghaderi and R~Srikant.
\newblock Opinion dynamics in social networks: A local interaction game with
  stubborn agents.
\newblock In {\em 2013 American Control Conference}, pages 1982--1987. IEEE,
  2013.

\bibitem{mirtabatabaei2012opinion}
Anahita Mirtabatabaei and Francesco Bullo.
\newblock Opinion dynamics in heterogeneous networks: convergence conjectures
  and theorems.
\newblock {\em SIAM Journal on Control and Optimization}, 50(5):2763--2785,
  2012.

\bibitem{williams1991probability}
David Williams.
\newblock {\em Probability with martingales}.
\newblock Cambridge university press, 1991.

\end{thebibliography}
\bibliographystyle{unsrt}}

\end{document}